\documentclass[12pt]{asl}

\begin{document}

\theoremstyle{plain}
\newtheorem{thrm}{Theorem}
\newtheorem{lmm}[thrm]{Lemma}
\newtheorem{crllry}[thrm]{Corollary}
\newtheorem{cnjctr}[thrm]{Conjecture}
\newtheorem{qstn}[thrm]{Question}

\theoremstyle{definition}
\newtheorem{dfntn}{Definition}

\newcommand{\ignore}[1]{}
\newcommand{\fld}[1]{\mathbb{#1}}

\author{Tobias J. Hagge}
\revauthor{Hagge, Tobias J.}

\address{Department of Mathematics, Indiana University, Bloomington, IN 47405}

\email{thagge@indiana.edu}

\ignore{\date{\today}}

\title{$QL(\fld{C}^n)$ determines $n$}

\begin{abstract}
This addendum to \cite{dhmw} shows that the set of tautological quantum logical propositional formulas for a finite dimensional vector space $\fld{C}^n$ is different for every $n$, affirmatively answering a question posed therein.
\end{abstract}

\maketitle
The paper \cite{dhmw} explored the properties of Birkhoff and Von Neumann's propositional quantum logic (see \cite{bvn}) as modelled by finite dimensional Hilbert spaces. One question asked in \cite{dhmw} is whether the set of tautological propositional formulas uniquely determines the dimension of the underlying vector space. A partial answer was given, namely that $\fld{C}^n$ and $\fld{C}^{2n}$ give different sets of tautologies. This note gives a full answer to the question.

For our purposes, propositional formulas consist of alphabet symbols, perentheses, and the symbols meet ($\land$), join ($\lor$), orthocomplement ($\lnot$), top ($\top$), and bottom ($\bot$). The well formed formulas are the same as those of propositional boolean logic. The symbol $\top$ is interpreted as a finite dimensional Hilbert space, $\bot$ is the trivial subspace, alphabet symbols are variables standing for vector subspaces of $\top$, $\land$ is intersection, $\lor$ is span of union, and $\lnot$ is orthogonal complement in $\top$. With these operations, the set of subspaces of $\top$ forms a bounded modular ortholattice.

Let $\bar v = v_1, \ldots, v_k$ be a list of alphabet symbols and let $\bar S = S_1, \ldots, S_k$ be a collection of subspaces of a finite dimensional Hilbert space $U$. Given a well formed formula $\phi(\bar v)$, the valuation $\Xi_U(\phi(\bar v), \bar S)$ is the subspace resulting from instantiating each $v_i$ with the subspace $S_i$ and performing the operations described by $\phi$ with universal space $U$. As a shorthand the valuation may be implicit; for example if $S$ and $T$ are subspaces of $U$ then $\Xi_U(v \land w, S, T)$ is abbreviated $S \land T$, and $U$ is inferred from context.

\begin{dfntn}
Let $\phi(\bar v)$ be a well formed formula. Define $\bar d_\phi: \fld{N} \to \fld{N}$ such that
\[
\bar d_{\phi(\bar v)}(n) = \max_{\bar S}(dim(\Xi_{\fld{C}^n}(\phi(\bar v), \bar S))).
\]
\end{dfntn}

\begin{dfntn}
A well formed formula $\phi(\bar v)$ is a {\em tautology} in $\fld{C}^n$ if $\bar d_\phi(n) = 0$.
\end{dfntn}

\ignore{By duality under orthocomplement, one could instead define {\em lowest dimension} in the obvious way and say that $\phi$ is a tautology in $\fld{C}^n$ if the lowest dimension function evaluates to $n$. However, when using the restriction process defined below, the present definitions are more convenient.}

\begin{dfntn}
$QL(\fld{C}^n)$ is the set of tautologies when $\top = \fld{C}^n$.
\end{dfntn}

The goal is to establish the following:

\begin{thrm}
$m < n \Rightarrow QL(\fld{C}^n) \subsetneq QL(\fld{C}^m)$.
\end{thrm}

In \cite{dhmw} it was shown that $QL(\fld{C}^{n+1}) \subset QL(\fld{C}^n)$, and $QL(\fld{C}^n) \ne QL(\fld{C}^{2n})$. A formula $\phi_1$ was constructed such that $\bar d_{\phi_1} = \lfloor \frac{n}{2} \rfloor$. From $\phi_1$ formulas $\phi_k$ were constructed in stages such that $\bar d_{\phi_k}= \lfloor \frac{\bar d_{\phi_{k-1}}}{2} \rfloor$. Thus $\phi_k$ is a tautology in $QL(\fld{C}^n)$ iff $log_2(n) < k$. Each $\phi_k = \phi_1 \vert_{\phi_{k-1}}$, which is defined as follows:\footnote{This definition corrects a small mistake in the definition given in \cite{dhmw}.}

\begin{dfntn}
Suppose that $\alpha(\bar{u})$ is a formula in $k$ variables $\bar{u} = u_1, \ldots u_k$. For convenience, use De Morgan's laws to replace $\alpha$ with an equivalent formula, also called $\alpha$, in which all negations are negations of atomic variable symbols. Given a second formula $\beta(\bar{v})$,  let $\alpha(\bar{u}) \vert_{\beta(\bar{v})}$ denote the modification of $\alpha$ such that each unnegated instance of $u_i$ is replaced with $u_i \land \beta(\bar{v})$ and each instance of $\lnot u_i$ is replaced with $\lnot (u_i \land \beta(\bar{v})) \land \beta(\bar{v})$.
\end{dfntn}

\begin{lmm}\label{lmmrstrct}
Let $\bar u = u_1, \ldots, u_{k_u}$ and $\bar v = v_1, \ldots, v_{k_v}$ be lists of variables, $\bar S = S_1, \ldots, S_{k_u}$ and $\bar T = T_1, \ldots, T_{k_v}$ lists of subspaces of $\fld{C}^n$, and $\alpha(\bar u)$ and $\beta(\bar v)$ formulas. Let $\bar P = P_1, \ldots P_{k_u}$ such that $P_i = S_i \land \Xi_{\fld{C}^n}(\beta(\bar v), \bar T)$. Then the following holds:
\[
\Xi_{\fld{C}^n}(\alpha(\bar u) \vert_{\beta(\bar v)},\bar S, \bar T) = \Xi_{\Xi_{\fld{C}^n}(\beta(\bar v), \bar T)}(\alpha(\bar u), \bar P)
\]
\end{lmm}
\begin{proof}
The construction procedure gives the result for atomic formulas and their negations. Since unions and intersections are not changed by inclusion into a larger universal space, the result follows by structural induction.
\end{proof}

\begin{crllry}\label{dmnsncrllry}
If $\bar d_{\alpha(\bar u)}=f$, and $\bar d_{\beta(\bar v)} = g$, then $\bar d_{\alpha(\bar u) \vert_{\beta(\bar v)}} = f \circ g$. In particular, $\alpha(\bar u) \vert_{\beta(\bar v)}$ is a tautology in $\fld{C}^n$ iff $\alpha(\bar u)$ is a tautology in $\fld{C}^{g(n)}$.
\end{crllry}

Because the function $n \to \lfloor \frac{n}{2} \rfloor$ is not injective, none of the formulas constructed in \cite{dhmw} distinguish dimensions between $2^k$ and $2^{k+1}-1$. To overcome this limitation, suppose $2^k \le m < n \le 2^{k+1}-1$ for some $k$, and assume there exists a formula $\alpha$ such that $\bar d_\alpha = \lfloor \frac{n}{2} \rfloor$. Construct a formula $\phi$ in stages, starting with $\phi_0 = \top$.

Suppose it is the beginning of stage $s$, and $\bar d_{\phi_{s-1}}(m) < \bar d_{\phi_{s-1}}(n)$. If $m$ is odd, or $\bar d_{\phi_{s-1}}(n)-\bar d_{\phi_{s-1}}(m) > 1$, define $\phi_s(\bar u, \bar v) = \alpha(\bar u) \vert_{\phi_{s-1}(\bar v)}$. Then by Corollary~\ref{dmnsncrllry}, $\bar d_{\phi_s} = \lfloor \frac{\bar d_{\phi_{s-1}}}{2} \rfloor$, and $\bar d_{\phi_s}(m) < \bar d_{\phi_s}(n)$. If, on the other hand, there is some $l$ such that $\bar d_{\phi_{s-1}}(m)=2l$ and $\bar d_{\phi_{s-1}}(n)=2l+1$, it will be shown that there is a formula $\beta_l$, which depends on $l$, such that $\bar d_{\beta_l}(2l)=l$ and $\bar d_{\beta_l}(2l+1)=l+1$. Then define $\phi_s(\bar u, \bar v) = \beta_l(\bar u) \vert_{\phi_{s-1}(\bar v)}$. Then by Corollary~\ref{dmnsncrllry},
\[
\bar d_{\phi_s}(m) = \bar d_{\beta_l}(\bar d_{\phi_{s-1}}(m)) = l < l+1 = \bar d_{\beta_l}(\bar d_{\phi_{s-1}}(n)) = \bar d_{\phi_s}(n).
\]

Since at each stage $\bar d_{\phi_s}(m) = \lfloor \frac{d_{\phi_{s-1}}(m)}{2} \rfloor$ and $\bar d_{\phi_s}(m) < \bar d_{\phi_s}(n)$, the construction procedure must eventually give $\phi_s$ such that $\bar d_{\phi_s}(m) = 0$ but $\bar d_{\phi_s}(n) > 0$. It remains only to construct $\alpha$ and $\beta_l$. A suitable formula for $\alpha$ was given in \cite{dhmw}, but here a new $\alpha$ will be constructed and then modified to give $\beta_l$.

Let $P_c$ denote the linear operator that projects onto the subspace given by the variable $c$. The following formula evaluates to the image of $P_b \circ P_a$:
\[
P(a,b) = (a \lor \lnot b) \land b.
\]

In a distributive lattice $P(a,b)=a \land b$, but in a modular lattice one only has $P(a,b) \ge a \land b$. The following are true when the lattice is subspaces of $\fld{C}^n$:

\begin{lmm} for all subspaces $S$ and $T$, the following hold:
  \begin{enumerate}
  \item $dim(P(S,T)) = dim(P(T,S)) \le min(dim(S), dim(T))$,
  \item $S \land P(S,T) = T \land P(T,S) = P(S,T) \land P(T,S) = S \land T$.
  \end{enumerate}
\end{lmm}

The proof is easy and omitted.

Define the formula
\[
\alpha(a,b) = P(b,a) \land \lnot (a \land b).
\]

When the distributive law holds $\alpha$ is a tautology, but $\alpha$ is not a tautology in $\fld{C}^n$ for $n \ge 2$.

\begin{lmm}\label{nvrtwlmm}
The formula $\alpha(a,b)$ has $\bar d_{\alpha} = \lfloor \frac{n}{2} \rfloor$ and for spaces $S$ and $T$, $dim(\Xi_{\fld{C}^n}(\alpha(a,b),S,T)) = \frac{n}{2}$ iff the following hold:
\begin{enumerate}
\item $n$ is even,
\item $dim(S)=dim(T)=\frac{n}{2}$,
\item $S \land T = \bot$,
\item $S \land \lnot T = \bot$.
\end{enumerate}

\end{lmm}
\begin{proof}
It is easy to verify that the above conditions on $S$ and $T$ give dimension $\frac{n}{2}$. For the other direction, if $min(dim(S),dim(T)) < \frac{n}{2}$, $dim(\alpha(S,T)) \le dim(P(S,T)) < \frac{n}{2}$. Also, since $S \land T \subset P(T,S)$, one gets the following: 

\begin{gather}
\label{smqtns} dim(\alpha(S,T)) = dim(P(T,S)) - dim(S \land T) \\
\le min(dim(T),dim(S)) - dim(S \land T) \\
\le min(dim(T),dim(S)) - dim(S) - dim(T) + dim(\top) \\
= dim(\top) - max(dim(S),dim(T)).
\end{gather}

Therefore, if $dim(\alpha(S,T))=\frac{n}{2}$, $dim(S)=dim(T)=\frac{n}{2}$. Thus $dim(P(T,S)) \le \frac{n}{2}$, so line~\ref{smqtns} implies that $dim(P(T,S)) = \frac{n}{2}$ and $dim(S \land T)=0$. Since $dim(P(T,S)) = \frac{n}{2}$, $dim(S \land \lnot T) = 0$.
\end{proof}

\begin{crllry}
$dim(\alpha(S,T)) = \frac{n}{2} \Rightarrow \alpha(S,T)=P(T,S)$.
\end{crllry}

To define $\beta_l$, restrict $\alpha$ to itself $\lfloor log_2(l) \rfloor-1$ times to obtain a formula $\gamma$ such that $\bar d_{\gamma}(2l) = \bar d_{\gamma}(2l+1)=1$. Define $\tilde \beta_l(a,b,\bar c) = \lnot (P(b,a) \lor P(a,b)) \land \gamma(\bar c)$ and $\beta_l(a,b,\bar c) = \tilde \beta_l(a,b,\bar c) \lor \alpha(a,b)$.

\begin{lmm}
The formula $\beta_l$ satisfies $\bar d_{\beta_l}(2l)=l$ and $\bar d_{\beta_l}(2l+1)=l+1$.
\end{lmm}
\begin{proof}
Clearly, $\bar d_{\tilde \beta_l}(2l) = \bar d_{\tilde \beta_l}(2l+1) = 1$. The conditions of Lemma~\ref{nvrtwlmm} imply that $\bar d_{\beta_l}(2l) = \bar d_{\alpha}(2l) = l$, while $\bar d_{\beta_l}(2l+1) = \bar d_{\alpha}(2l+1) + 1 = l+1$.
\end{proof}

\bibliographystyle{asl}
\bibliography{references}

\end{document}